\documentstyle[12pt]{report}

\renewcommand\a{\alpha}

\renewcommand\c{\circ}
\newcommand\C{{\mbox{\rm\bf C\hspace{-7.3pt}}{^{_{\bf\mid}}}\hspace{4.5pt}}}

\renewcommand\d{\delta}
\newcommand\dd{\Delta}
\newcommand\D{{\mbox{\rm\bf D\hspace{-7.3pt}}{^{_{\bf\mid}}}\hspace{4.5pt}}}

\newcommand\g{\gamma}

\renewcommand\ll{\lambda}

\newcommand\M{{\cal M}}

\newcommand\oo{\omega}

\newcommand\op[1]{\mathop{\rm #1}\nolimits}
\newcommand\p{\partial}
\newcommand\po{$\!\!\!{\bf .}$ }

\newcommand\R{{\rm I\hspace{-2.5pt} R}}
\renewcommand\t{\tau}
\newcommand\ti{\times}
\newcommand\te{\theta}
\newcommand\T{\Theta}

\newcommand\ve{\varepsilon}

\newcommand\V{{\cal V}}

\newcommand\x{\xi}

\newcommand\z{\zeta}
\newcommand\Z{{\rm Z\mkern-5muZ}}

\def\Rom#1{\uppercase\expandafter{\romannumeral#1}}

\newcommand\1{{\bf 1}}

\newcommand\qed{\phantom{\underline{y}}\hfill\hfill$\Box$}

\newcommand\bib[1]{\bibitem[#1]{#1}}

\newcommand{\text}[1]{{\mbox{\rm #1}}}
\newcommand{\dfrac}[2]{\frac{\displaystyle #1}{\displaystyle #2}}

\newtheorem{th}{Theorem}
\newtheorem{prop}{Proposition}

{\endtrivlist}
{\endtrivlist}
\newenvironment{dfn}{\trivlist \item[\hskip \labelsep{\bf Definition.}]}%
{\endtrivlist}
\newenvironment{lem}{\trivlist \item[\hskip \labelsep{{\bf Lemma.}}]}%
{\endtrivlist}
\newenvironment{proof}{\trivlist \item[\hskip
\labelsep{{\it\underline{Proof}.\/}}]}%
{\endtrivlist}
{\endtrivlist}

\newcounter{a}
\newcounter{f}
\makeatletter
\newcommand{\@thefnmark}{$^\fnsymbol{f}$}
\renewcommand{\@makefnmark}{\hbox{\mathsurround=0pt
                           $^{\fnsymbol{f}}$}}
\renewcommand{\@makefntext}[1]{\parindent=1em\noindent
            \hbox to 1.8em{\hss$^{\fnsymbol{f}}$}#1}
\makeatother

\begin{document}

\title{On the Kobayashi-Royden pseudonorm \\
       for almost complex manifolds}

\author{\bf Boris~S.~Kruglikov}

\date{August 15, 1997}

\maketitle

\begin{abstract}
In this paper we define Kobayashi-Royden pseudonorm for almost complex
manifolds. Its basic properties known from the complex analysis are
preserved in the nonintegrable case as well. We prove that the
pseudodistance induced by this pseudonorm coincides with the Kobayashi
pseudodistance defined for the almost complex case earlier. We also consider
a geometric application for moduli spaces of pseudoholomorphic curves.
\end{abstract}

\tableofcontents

\clearpage

%%%%%%%%%%%%%%%%%%%%%%%%%%%%%%%%%%%%%%%%%%%%%%%%%%%%%%%%%%%%%%%%%%%%%%%%%%

\chapter*{Introduction}
\addcontentsline{toc}{chapter}{\bf\quad \  Introduction}

\hspace{13.5pt}
In 1967 Kobayashi (\cite{K1}) introduced biholomorphically-invariant
pseudodistance (which is a distance without nondegenerace axiom)
for any complex manifold. This gave rise for the theory of hyperbolic spaces
(see~\cite{K2,La,PS} for details). The Kobayashi pseudodistance is the
maximal pseudodistance among all nonincreasing under holomorphic maps
pseudodistances and such that on the unit disk $\D\subset\C$ it coincides
with the distance $d_\D$ which is induced by infinitesimal Poincar\'e metric
of the Lobachevskii plane:
 $$
dl^2=\frac{dz d\bar z}{(1-|z|^2)^2}.
 $$
Another way to define the pseudodistance $d_M$ on a complex manifold
$M$ is the following:
 $$
d_M(p,q) = \op{inf} \ \sum_{k=1}^{m}d_\D(z_k,w_k),
 $$
where the infimum is taken over all holomorphic mappings
$f_k: \D\rightarrow M$, $k=1,...,m$, with the properties
$f_1(z_1) = p$, $f_{k}(w_k) = f_{k+1}(z_{k+1})$ and $f_m(w_m) = q$.
In the paper~\cite{KO} the Kobayashi pseudodistance notion was extended
to the case of arbitrary almost complex manifolds and it was shown that 
basic properties of this pseudodistance stay preserved.

In 1970 Royden (\cite{Ro}) discovered and justified an infinitesimal
analog of the Kobayashi pseudodistance for complex manifolds. The main goal
of the present paper is to define the corresponding notion in the category
of almost complex manifolds and to prove the coincidence theorem. We also
consider the reduction procedure which permits to define geometrical
invariants for moduli spaces of pseudoholomorphic curves.%
\footnote{
The author thanks all scientists from Math. and Stat. Dept. of 
the Univ. of Tromsoe where he was a visitor and where the paper 
was finished.}

%%%%%%%%%%%%%%%%%%%%%%%%%%%%%%%%%%%%%%%%%%%%%%%%%%%%%%%%%%%%%%%%%%%%%%%%%%
% 1 %
 \chapter{Definition of the pseudonorm and basic properties}

\hspace{13.5pt}
Let $(M^{2n},J)$ be an almost complex manifold, i.e. $J^2=-\1\in
T^*M\otimes TM$. We denote by $e=1\in T_0\D$ the unit vector.
Let us also denote by $\D_R$ the disc of radius $R$ in $\C$. Let
$v\in T_pM$, $p=\t_Mv$, where by $\t_M: TM\to M$ we denote the canonical
projection.

 \begin{dfn}
The Kobayashi-Royden pseudonorm is the function
 $$
F_M(v) = \inf\limits_{{\cal{R}}}\frac{1}{r},
 $$
where the set ${\cal{R}}$ consists of all pseudoholomorphic mappings
$f:\D\to M$ (i.e. such mappings that $f_*\c j_0=J\c f_*$) with the
property $f_*(0)e=rv$, $r\in\R_+$.
 \end{dfn}

By theorem III from~\cite{NW} the set ${\cal{R}}=\bigcup{\cal{R}}_r(v)$
is nonempty since ${\cal{R}}_r(v)$ is nonempty for every $r>0$ small enough.
Here ${\cal{R}}_r(v)$ stands for the set of all pseudoholomorphic curves $f$
with $f_*(0)e=rv$, $r>0$.
Next statement is the direct consequence of the definition.

 \begin{prop}\po
For any vector $v\in T_*(M_1)$ and for any pseudoholomorphic mapping
$f: (M_1,J_1)\to (M_2,J_2)$ we have:
 $$
\phantom{aaaaabbbbbbb}\qquad\hfill\qquad\qquad\hfill
F_{M_2}(f_*v)\le F_{M_1}(v).
\phantom{aaaaabbbbbbb}\qquad\hfill\qquad\qquad\hfill\Box
 $$
 \end{prop}

 \begin{prop}\po
{\rm(i)}. The function $F_M$ on $TM$ is nonnegative and homogeneous
of degree one: $F_M(tv)=|t|F_M(v)$ for any $v\in T_*M$ and $t\in\R$. \newline
{\rm(ii)}. Let $K\subset M$ be a compact with nonempty set of interior
points such that the almost complex structure $J$ is tamed by an exact
symplectic form $\oo=d\a$ (i.e. $\oo(\x,J\x)>0$ for $\x\ne0$) on it.
Then for any $v\in T_pM$ with $p=\t_M v\in K$ and for any norm $|\cdot|$ on
$\tau_M$ there exists a constant $C_K$ such that
 $$
F_M(v)\le C_K |v|.
 $$
 \end{prop}

 \begin{proof}
The first statement is obvious while the second is a reformulation of the
nonlinear Schwarz lemma (\cite{Gr1}). \qed
 \end{proof}

 \begin{prop}\po
The function $F_M$ is upper semicontinuous.
 \end{prop}

 \begin{proof}
Our statement is equivalent to the following one: if there is a
pseudoholomorphic disk $f: \D_R\to M$ in some direction $v\in T_*M$,
$f_*e=v$, then in any close direction $v'\in\V(v)$ of $v$ there is a
pseudoholomorphic disk of almost the same size $f': \D_{R-\ve}\to M$,
$(f')_*e=v'$. Here $\ve>0$ is arbitrary small fixed number. The last
statement is equivalent to the existence theorem for some nonlinear partial
differential equation. We prove it using the fixed point theorem in
the Banach space following~\cite{NW}. This paper by theorem III provides
existence of a small pseudoholomorphic disk in any given direction.
The proof is based on the linearization of the almost complex structure at
the point, on the writing down the corresponding nonlinear equation and on
the proximity of the obtained equation to the linear equation on holomorphic
curves of the linearized structure. In our case we can linearize the almost
complex structure along given pseudoholomorphic curve $f: \D_R\to M$.

Contracting if needed we can assume that our manifold is a neighborhood of
the disk $\D_R\subset\C\ti\{0\}^{n-1}\subset\C^n$ and that
the almost complex structure $J$ coincides  with the standard complex
structure $J_0$ at the points of this disk. Let $k\in\Z_+$, $\ll\in(0,1)$.
Denote by $C^{k+\ll}_{(R)} = C^{k+\ll}(\D_R;\C^n)$ the space of all
$\ll$-H\"older $k$-smooth mappings of $\D_R$ into $\C^n$ and consider
this space equipped with the
standard H\"older norm $|\!|\cdot|\!|$. Let's introduce also pseudonorm
$|\!|f|\!|'=\max\{ |\!|\p f|\!|,\ |\!|\bar\p f|\!|\}$.
We use standard notations $\p_i$ and ${\bar\p}_i$ defined regardness the
complex structure $J_0$. In this case the defining equation for
pseudoholomorphic curve $f: \D_{R-\ve}(\z)\to \C^n(z^i)$ in coordinates
$z^i=z^i(\z)$ takes the form (see details in~\cite{NW}):

 $$
\bar\p z^i+ \sum_m a^i_{\bar m}(z)\bar\p {\bar {z}}^m =0,
\qquad z^i\in C^{k+\ll}(\D_{R-\ve};\C^n). \eqno (*)
 $$
Fixation of the point $p\in\C^n$ which the curve goes through and
of the direction of the curve $v\in T_pM$ takes the form of initial
conditions: $\ z^i(0)=p^i$, $\p z^i(0)=v^i$. Due to our linearization
we have $a^i_{\bar m}(z)=0$ for all points $z\in \D_R\subset \C^n$.
Thus we have the solution $z_0=\z v_0$ in the direction of the vector
$v_0=(1,0,\dots,0)\in T_0\C^n$, $|\z|<R$. We seek for a solution of the
equation $(*)$ at the form
 $$
z(\z) = p + \z v + \T(z,z),
 $$
 $$
\T^i(f,g)(\z)=\te^i(f,g)(\z)-\te^i(f,g)(0)-\z[\p\te^i(f,g)](0),
 $$
 $$
\te^i(f,g)=-T\left(\sum_m a^i_{\bar m}(g) \bar\p {\bar f}^m\right),
 $$
where $Tf(w)=\dfrac1{\mathstrut 2\pi i}{\displaystyle\int_{D_R}}
\dfrac{f(\z)}{\z-w}d\z\wedge d\bar\z$.
Let us suppose that $a^i_{\bar m}\in C^{k+\ll}$. Since the function
$z_0$ is a solution of our equation then $\T(z_0,z_0)=0$. Consider
its neighborhood ${\cal{O}}\subset C^{k+\ll}_{(R-\ve)}$. Consider also
a neighborhood ${\cal{V}}$ of the vector $v_0$ in $T_*\C^n$.
For sufficiently small neighborhood ${\cal{V}}$ we seek a solution
of the equation $(*)$ in ${\cal{O}}$ as the limit of the iteration process

 $$
z^{(v)}_1= p +\z v,\ \
z^{(v)}_{k+1}= p +\z v + \T(z^{(v)}_k,z^{(v)}_k).
 $$

 \begin{lem}
For any given small $\ve>0$ the iterative process converges in
$C^{k+\ll}_{(R-\ve)}$ for sufficiently small neighborhood
${\cal{V}}$ and the limit is a function $z^{(v)}\in C^{k+1+\ll}_{(R-\ve)}$.
This function is a solution of the equation $(*)$ with initial conditions
$z^{(v)}(0)=p$, $\p z^{(v)}(0)=v$, and in addition $z^{(v)}\to z_0$ when
$v\to v_0$ and for small enough ${\cal{V}}$ we have:
$|\!|z^{(v)}|\!|'\le 2|v|$.
 \end{lem}

The required convergence is proved by estimates similar to~\cite{NW} but in
our case the polidisk has the form
$\D_{R-\ve}\ti(\D_{\delta(\ve)})^{n-1}\subset \C^n$.
The main estimate is the upper bound of the contraction coefficient of the
process:
 $$
|\!|\T(f,f)-\T(g,g)|\!|<c\cdot|v-v_0|^\ll\cdot|\!|f-g|\!|.
 $$
It is obtained via the following estimates (cf.~\cite{NW}):

 \begin{enumerate}
  \item
$|\!|Th|\!|'\le c_1\cdot|\!|h|\!|$.
  \item
$|\!|\T(f,f)-\T(g,g)|\!|\le
|\!|\T(f,f)-\T(f,g)|\!| + |\!|\T(f,g)-\T(g,g)|\!|$.
  \item
$|\!|a^i_{\bar m}(f)|\!|\le c_2\cdot|v-v_0|^\ll$.
  \item
$|\!|a^i_{\bar m}(f)-a^i_{\bar m}(g)|\!|\le c_3\cdot|v-v_0|\cdot
|\!|f-g|\!|$.
 \end{enumerate}
$\phantom{aaaaa}$ \qed
 \end{proof}

%%%%%%%%%%%%%%%%%%%%%%%%%%%%%%%%%%%%%%%%%%%%%%%%%%%%%%%%%%%%%%%%%%%%%%%%%%
% 2 %
 \chapter{Coincidence theorem}

\hspace{13.5pt}
Let us define the function ${\bar d}_M: M\ti M\to \R$ by the formula
 $$
{\bar d}_M(p,q)=\inf\limits_\g \int_0^1 F_M(\dot\g(t))\,dt
 $$
where the infimum is taken over all piecewise smooth paths from the point
$p$ to the point $q$. The upper semicontinuity and upper boundness
(propositions 2 and 3) gives us the correctness of this definition and

 \begin{prop}\po
The function ${\bar d}_M$ is a pseudodistance. \qed
 \end{prop}

 \begin{th}\po
The pseudodistance just introduced coincides with the Kobayashi
pseudodistance: $d_M={\bar d}_M$.
 \end{th}

 \begin{proof}
The inequality ${\bar d}_M\le d_M$ is obvious because $F_M(v)=\inf |\xi|$
where the infimum is taken over all pseudoholomorphic mappings
$f: \D\to M$, $f_*\xi=v$, and the norm is induced by the Poincar\'e metric.
Let's prove the inverse. We follow the Royden's proof (\cite{Ro}).

Let $\g$ be such a smooth path from a point $p$ to a point $q$ that
${\displaystyle\int_\g}F_M< {\bar d}_M(p,q)+\ve$. Due to upper
semicontinuity there exists a continuous function $h$ on the interval
$[0,1]$ such that $h(t)>F_M(\dot \g (t))$ and
 $$
\int_0^1 h(t)\,dt < {\bar d}_M(p,q)+\ve,
 $$
i.e. for a sufficiently dense partition $0=t_0<t_1<\dots <t_k=1$ we have:
 $$
\sum_{i=1}^k h(t_{i-1}) (t_i-t_{i-1})< {\bar d}_M(p,q)+\ve.
 $$

Let us consider arbitrary pseudoholomorphic curve $u_t^\g: \D_\d\to M$
which satisfies the conditions $u_t^\g(0)=\g(t)$ and
$(u_t^\g)_*e=\dot\g (t)$. Define for small $\dd t\in\R_+\subset\C$ the curve
$\hat\g(t;\dd t)=u_t^\g(\dd t)$. Since
$\hat\g (t;\dd t)=\g(t+\dd t)+ O(|\dd t|^2)$ then due to proposition~2 for
small enough $\dd t$ we have:
  \begin{eqnarray*}
d_M(\g(t),\g(t+\dd t))
&\!\!\!\le\!\!\!&
d_M(\g(t),\hat\g(t;\dd t)) + d_M(\hat\g(t;\dd t), \g(t+\dd t)) \\
&\!\!\!\le\!\!\!&
F_M(\dot\g(t))\dd t + O(|\dd t|^2)\le (1+\ve)h(t)\dd t.
  \end{eqnarray*}

Thus for sufficiently dense partition of the interval $[0,1]$ we have:
  \begin{eqnarray*}
d_M(p,q)
&\!\!\!\le\!\!\!&
\sum_{i=1}^k d_M(\g(t_{i-1}),\g(t_i)) \\
&\!\!\!\le\!\!\!&
(1+\ve)\sum_{i=1}^k h(t_{i-1})(t_i-t_{i-1})<
(1+\ve)({\bar d}_M(p,q)+\ve).
  \end{eqnarray*}
Since $\ve>0$ is arbitrary we are done. \qed
 \end{proof}

%%%%%%%%%%%%%%%%%%%%%%%%%%%%%%%%%%%%%%%%%%%%%%%%%%%%%%%%%%%%%%%%%%%%%%%%%%
% 3 %
 \chapter{Connection with $h$-principle}

\hspace{13.5pt}
Let us call a pseudoholomorphic disks fiber bundle any embedding
(immersion) $\Phi: \D_R\ti N^{2n-2}\to M$ such that all the mappings
$\left.\Phi\right|_{D_R\ti\{x\}}$ are pseudoholomorphic.

 \begin{prop}\po
Let $(M,J)$ be a 4-dimensional almost complex manifold. For any embedded
(immersed) pseudoholomorphic disk $f: \D_R\to M$ and any $\ve>0$ a small
neighborhood of the image $f(\D_{R-\ve})$ can be fibered by
pseudoholomorphic disks.
 \end{prop}

 \begin{proof}
We can change the almost complex structure $J$ in a neighborhood of the
boundary of the disk $f(\p\D_R)$ turning it into the standard complex.
After this we can reduce the manifold $M$ gluing up the disk into the
sphere. Our manifold takes the form $S^2_R\ti N^2$ (the subscript
for the sphere has the following meaning: for our future purposes we
introduce a symplectic form $\oo$ taming $J$, $\oo(\x,J\x)>0$, of the form
$\oo=\oo_1\oplus\oo_2$, the volume being $\op{vol}_{\oo_1}(S^2_R)=\pi R^2$).
We may also change the almost complex structure everywhere outside a small
neighborhood of the disk image turning it into the standard integrable.
After the corresponding gluing up we may assume that our manifold takes the
form $M^4=S^2_R\ti S^2_{2R}$ and in addition the almost complex structure
$J$ is the standard complex out of a small neighborhood of
$f(\D_{R-\ve})$. Denote by $A\in H_2(M;\Z)$ the homological class of
the sphere $S^2_R\ti\{*\}\subset M$. According to~\cite{Gr1} for an almost
complex structure $J$ of a general position there exists a unique
pseudoholomorphic $A$-sphere through any point of the manifold $M$.
These spheres form in totality a pseudoholomorphic disks fiber bundle.
The absence of intersections follows from nonnegativity of the intersection
number (\cite{Gr1}), this construction is similar to lemmas
3.5 and 4.1 from~\cite{M1} where the smoothness of the fibration is also
proved. The intersection of the obtained bundle with a small neighborhood
of the image $f(\D_{R-\ve})$ gives  us a fibration of this image
(with the initial almost complex structure $J$) by pseudoholomorphic disks.
In order to get rid of the general position condition for the structure
$J$ on $S^2_R\ti S^2_{2R}$ we perturb it and then go to limit
with the perturbation tending to zero and use the compactness theorem
(\cite{MS}). \qed
 \end{proof}

According to~\cite{Gau} (cf.~\cite{Kr}) for any almost complex
$(M^{2n},J)$ there exists a  canonical almost complex structure $J_{[1]}$
on the corresponding manifold of pseudoholomorphic 1-jets $J^1_{PH}(\D_R;M)$
such that the canonical projection $\pi: J^1_{PH}(\D_R;M)\to M$ is a
pseudoholomorphic mapping. Any pseudoholomorphic mapping $f: \D_R\to M$
can be canonically lifted to a pseudoholomorphic mapping
$j^1f: \D_R\to J^1_{PH}(\D_R;M)$.

In any even dimension $2n$ similarly to proposition~5 one can prove that for
any embedded (immersed) pseudoholomorphic disc $f: \D_R\to M^{2n}$ and any
sufficiently close to its image point there exists a pseudoholomorphic disk
$f': \D_{R-\ve}\to M^{2n}$ through this point. We may consider
such a map $g': \D_{R-\ve}\to J^1_{PH}(\D_{R-\ve};M)$ constructed by the
map $g=j^1f$. Note that the constructed map is not necessary the jet lift
of some map, $g'\ne j^1(\pi\circ g')$. If the equality $g'=j^1f'$ holds
then, in accordance with the terminology from~\cite{Gr2}, the corresponding
map should be called {\it holonomic\/}. Closeness of the initial points of
the disks $g$ and $g'$ (i.e. the points $g(0)$ and $g'(0)$) in
$J^1_{PH}(\D_{R-\ve};M)$ means closeness of initial points and
initial directions for the maps $f$ and $f'$ in $T_*M$. Thus
proposition~3 implies

 \begin{th}\po
There exists an embedded (immersed) holonomic pseudoholomorphic disk
$g': \D_{R-\ve}\to J^1_{PH}(\D_{R-\ve};M)$ through every point sufficiently
close to the image of an embedded (immersed) holonomic pseudoholomorphic
disk $g: \D_R\to J^1_{PH}(\D_R;M)$. \qed
 \end{th}

We may assume we have a family of pseudoholomorphic disks
$g_\a: \D_{R-\ve}\to J^1_{PH}(\D_{R-\ve};M)$ which fill out a neighborhood
of the image of $g_0=j^1f$. Under which conditions this family may be
assumed fibering (as in proposition~5)?

%%%%%%%%%%%%%%%%%%%%%%%%%%%%%%%%%%%%%%%%%%%%%%%%%%%%%%%%%%%%%%%%%%%%%%%%%%
% 4 %
 \chapter{Hyperbolicity and nonhyperbolicity}

\hspace{13.5pt}
Almost complex manifold $(M,J)$ is called ${\it hyperbolic\/}$ if the
pseudodistance $d_M$ is a distance. Let us consider arbitrary norm
$|\cdot|$ on the tangent bundle $\tau_M: TM\to M$ and let us consider the
corresponding bundle of unit tangent vectors $\tau^{(1)}_M: T_1M\to M$.
Consider the restriction of the Kobayashi-Royden pseudonorm on this
submanifold, $F_M^{(1)}: T_1M\to \R$. Proposition~2(ii) and
theorem~2(i) from~\cite{KO} imply

 \begin{th}\po
The function $F_M^{(1)}$ is bounded on compact subsets in $M$.
The manifold $M$ is hyperbolic if and only if $F_M^{(1)}\ne0$ and the
function $(F_M^{(1)})^{-1}$ is bounded on compact subsets in $M$, i.e.
$F_M^{(1)}$ is bounded away from zero on compact subsets. \qed
 \end{th}

Note that theorems 1 and 3 imply theorem~2(ii) from~\cite{KO}, that is the
almost complex version of Brody's theorem which states that for compact
(possible with boundary) manifold $M$ hyperbolicity is equivalent to the
absence of entire pseudoholomorphic curves, i.e. nontrivial
pseudoholomorphic mappings $f:\C\to M$.

Now let us consider the case of nonhyperbolic manifold $M$, for example let
it possess pseudoholomorphic spheres. In the case of general position for
the almost complex structure $J$, which is tamed by some symplectic form
$\oo$ on $M$, the set of all pseudoholomorphic spheres in a fixed homology
class $A\in H_2(M;\Z)$ (completed for compactness by the set of
decomposable rational curves) is a finite-dimensional manifold
$\M(A;J)$ (\cite{Gr1,MS}). We define by the {\it reduction procedure\/}
some pseudodistance on this manifold. Namely for any two
pseudoholomorphic spheres $f_i: S^2\to M$ defined up to holomorphical
reparametrization of $S^2$ let
 $$
d_{\M}([f_1],[f_2])=d_M(p_1,p_2), \eqno(\dagger)
 $$
where $p_i\in\op{Im}(f_i)$ are arbitrary points on the images.

As example note that defined pseudodistance is a distance for the manifold
$M^4=\Sigma^2_g\ti S^2$ with $g>1$ with some almost complex structure. 
Actually by methods of~\cite{Gr1}, as in proposition~5, one may prove that 
$M$ is fibered by pseudoholomorphic spheres and $\M\simeq \Sigma^2_g$.

Note that it is shown in the paper~\cite{M2} that for $N$ large enough
the manifold $\M\ti\R^{2N}$ possesses a homotopically canonical almost 
complex structure $\tilde J$. Applying the reduction procedure
we may assume that the Ko\-ba\-ya\-shi pseudodistance $d_{\M\ti\R^{2N}}$
induces a pseudodistance on $\M$. There is an existence question for almost
complex structures $\tilde J$ such that the induced pseudodistance on $\M$
equals to the reduced pseudodistance $d_\M$ defined by the formula 
$(\dagger)$.

%%%%%%%%%%%%%%%%%%%%%%%%%%%%%%%%%%%%%%%%%%%%%%%%%%%%%%%%%%%%%%%%%%%%%%%%%%%%
%
%%%%%%%%%%%%%%%%%%%%%%%%%%%%%%%%%%%%%%%%%%%%%%%%%%%%%%%%%%%%%%%%%%%%%%%%%%%%

\ {\hbox to 12.5cm{ \hrulefill }}

\bigskip

{\it Address:}
{\footnotesize
 \begin{itemize}
  \item
P.O. Box 546, 119618, Moscow, Russia
  \item
Chair of System Analysis, Moscow State Technological University
\linebreak
{\rm n. a.} Baumann, Moscow, Russia
 \end{itemize}
}

{\it E-mail:} \quad
{\footnotesize
lychagin\verb"@"glas.apc.org or borkru\verb"@"difgeo.math.msu.su
}

%%%%%%%%%%%%%%%%%%%%%%%%%%%%%%%%%%%%%%%%%%%%%%%%%%%%%%%%%%%%%%%%%%%%%%%%%%%%
\end{document}